\documentclass[12pt]{article}
\usepackage{amsmath}
\usepackage[light,firsttwo,bottom]{draftcopy}
\usepackage{amssymb}
\newcommand{\nc}{\newcommand}
\nc{\nt}{\newtheorem}
\nt{Satz}       {Theorem} [section]
\nt{Cor}    [Satz]  {Corollary}
\nt{Prop}   [Satz]  {Proposition}
\nt{Vermutung}  [Satz]  {Vermutung}
\nt{Bem}        [Satz]  {Remark}
\nt{Bs}     [Satz]  {Beispiel}
\nt{Lem}    [Satz]  {Lemma}       
\nt{Def}        [Satz]  {Definition}  
\nt{Maintheorem}        {Main Theorem}
\newcommand{\ultra}[4]
{
p_{#1,\underline{#2}}\left(#3,#4\right)
}
\newcommand{\diff}[2]
{
\cd_{#1,\underline{#2}}
}
\nc{\diffknu}{\diff {k}{2 \nu}}
\nc{\dimens}{\mathsf{dim}}
\newcommand{\textdiffjm}[2]
{$
\cd_{#1,\underline{#2}}^{\mathsf{J},m}
$}
\newcommand{\diffjm}[2]
{
\cd_{#1,\underline{#2}}^{\mathsf{J},m}
}
\nc{\gj}{\Gamma^J}
\nc{\buq}{``}
\nc{\eoq}{''}
\nc{\boq}{``}
\nc{\akzent}{\`}

\nc{\ca}{{\cal A}}
\nc{\cb}{{\cal B}}
\nc{\cc}{{\cal C}}
\nc{\cd}{{\cal D}}
 \nc{\cn}{{\cal N}}
\nc{\ci}{{\cal I}} \nc{\crr}{{\cal R}} \nc{\cj}{{\cal J}}
\nc{\kt}{K[T_1,\ldots,T_n]} 
\nc{\cf}{\cal F} \nc{\ch}
{\mbox{$\cal H$}} 
\nc{\kG}{\tilde{k}_G} \nc{\kg}{\tilde{k}_g}
\nc{\kh}{\tilde{k}_h} \nc{\kf}{\tilde{k}_f}\nc{\kF}{\tilde{k}_F}
 \nc{\td}{\tilde{\Delta}}
 \nc{\tildetau}{\tilde{\tau}}
 
\nc{\tk}{{\tilde{k}}}
\nc{\hutvontk}{{\widehat{d_{{\tk}}}}}
\nc{\hutvontknu}{{\widehat{d}_{{k + 2 \nu}}}}
\nc{\gaussk}{\gauss{\frac{k}{10}}}
\nc{\maass}{M_k^{\mathsf{Maass}}}
\nc{\smaass}{S_k^{\mathsf{Maass}}}
\nc{\mkzwei}{M_k^{(2)}}
\nc{\skzwei}{S_k^{(2)}}
\nc{\mk}{M_k}
\nc{\sk}{S_k}
\nc{\mknu}{M_{k+ 2\nu}}
\nc{\sknu}{S_{k + 2 \nu}}
\nc{\mksym}{M_k^{\mathsf{Sym}}}
\nc{\sksym} {S_k^{\mathsf{Sym}}}
\nc{\mksymnu}{M_{k + 2 \nu}^{\mathsf{Sym}}}
\nc{\mksymnunull}{M_{k + 2 \nu_0}^{\mathsf{Sym}}}
\nc{\sksymnu} {S_{k + 2 \nu}^{\mathsf{Sym}}}

\nc{\mkd}{M_k^{D}}
\nc{\skd}{S_k^{D}} 
\nc{\mkdnu}{M_{k + 2 \nu}^{D}}
\nc{\skdnu}{S_{k+ 2 \nu}^{D}}

\nc{\Mkzwei}{M_k^{(2)}}
\nc{\Skzwei}{S_k^{(2)}}
\nc{\kzwoelf}{\lfloor \frac{k}{12}\rfloor}
\nc{\knuth}[1]{\lfloor #1\rfloor}
\nc{\sedi}{\rtimes}
\nc{\ep}{\varepsilon}
\nc{\minuss}{\mbox{\scriptsize \mbox{-}}}      
\nc{\minus}{\mbox{\small \mbox{-}}}            
\nc{\minusv}[1]{\minus #1}
\nc{\p}{\ast}
\nc{\fein}{P_{4,0}(Q)\backslash Sp_4(Q)/Sp_1(Q)^4}
\nc{\ts}{\tau'}
\nc{\ys}{y_{\ts}}                              
\nc{\tz}{\left( \begin{array}{cc} \tau_1& 0\\ 0 &\tau_2 \end{array} \right)
                ,(z_1 z_2)}
\nc{\gu}{\Gamma_{\infty}}
\nc{\A}{\Bbb{A}}
\nc{\E}{\Bbb{E}}
\nc{\LL}{\Bbb{L}}
\nc{\HH}{\Bbb{H}}
\nc{\N}{\Bbb{N}}                         
\nc{\R}{\Bbb{R}}                     
\nc{\C}{\Bbb{C}}                      
\nc{\Q}{\Bbb{Q}}                     
\nc{\Z}{\Bbb{Z}}
\nc{\D}{\Bbb{D}}
\nc{\Hn}{I\!\!H_n}  \nc{\Ha}{I\!\!H}
\newcommand{\vier}[4]
{\left(\begin{array}{cc}#1 &#2 \\#3&#4 \end{array}\right)}
\newcommand{\ovier}[4]
{{\small\begin{array}{cc}#1 &#2 \\#3&#4 \end{array}}}

\newcommand{\tovier}[4]
{{\tiny\begin{array}{cc}#1 &#2 \\#3&#4 \end{array}}}

\nc{\einslambda}{[ \binom{1}{\lambda} ] }
\nc{\einslambdatext}{[ \tbinom{1}{\lambda} ] }
\nc{\einst}{[ \dbinom{1}{t} ] }
\nc{\einsttext}{[ \tbinom{1}{t} ] }
\newenvironment{beweis}{ {\bf Proof:}
\hfil\break}{\nopagebreak\hspace*{\fill}\raisebox{-2ex}
{\fbox{\rule{0mm}{0.6ex}\hspace{0.6ex}}}\\[4ex]}
\nc{\phiF}{\tilde{\Phi}_1^F}
\nc{\phiGG}{\tilde{\Phi}_{nl^2}^G}
\nc{\phiG}{\tilde{\Phi}_1^G}
\nc{\phiGpl}{\tilde{\Phi}_{p^{n+2l}}^G}
\nc{\phiGp}{\tilde{\Phi}_{p^n}^G}
\nc{\lk}{\langle}
\nc{\rk}{\rangle}
\nc{\llk}{{\Huge \langle}}
\nc{\rrk}{{\Huge \rangle}}
\nc{\kleinergleich}{\leqslant}
\nc{\grossergleich}{\geqslant}
\nc{\halbganz}{{\A}_2^{m}}
\nc{\mrn}{[n,r,m]}
\nc{\jkm}{J_{k,m}}
\nc{\jkmcusp}{J_{k,m}^{\, \mathrm{cusp}}}
\nc{\jknumcusp}{J_{k + 2 \nu,m}^{\, \mathrm{cusp}}}
\nc{\jkmEis}{J_{k,m}^{\, \mathrm{Eis}}}
\nc{\cdcusp}{\cd^{\, \mathrm{cusp}}}
\nc{\cdm}{\cd_{\nu }^{m}}
\nc{\cdnull}{\cd_{0}^{m}}
\nc{\cdzweim}{\cd_{2 \nu }^{m}}
\nc{\pkT}{P_{k,T}^J}
\nc{\Teinslambda}{T[\binom{1}{\lambda}]}
\nc{\Gegenbauernu}{{\bf \Huge \mathcal{C}_{2 \nu}^{k-1}}} 
\nc{\Gegenbauernull}{{\bf \Huge \mathcal{C}_{0}^{k-1}}}
\nc{\CC}{\mathcal{C}_{2 \nu}^{k-1}}
\nc{\CCjlambda}{\mathcal{C}_{2 j}^{\lambda}} 
\nc{\Satoh}{\{    \}_{k,m}^{\nu}}
\nc{\Satohpurdelta}{\{ \Delta  \}}
\nc{\Satohnull}{\{  \: \:  \}_{k,m}^{0}}
\nc{\Satohnix}{\{  \: \:  \}_{k,m}}
\nc{\Satohg}{\{  g  \}_{k,m}^{\nu}} 
\nc{\gauss}[1] {\lfloor #1 \rfloor} 
\begin{document}
\setlength{\parindent}{0cm}
\sf                     
\begin{center}{\Large \bf On The Spezialschar Of Maass           }
       \end{center}
    \begin{center}{By {\sc Bernhard Heim}}\\ 
Version 2.0: 16.02.2006
       \end{center}
\begin{abstract}
Let $M_k^{(n)}$ be the space of Siegel modular forms of degree $n$ and even weight $k$. 
In this paper firstly a certain subspace $\mathsf{Spez}(M_k^{(2n)})$ the 
Spezialschar of $M_k^{(2n)}$ is introduced. In the setting of the Siegel three-fold it is proven that this 
Spezialschar is the Maass Spezialschar. Secondly an embedding of $M_k^{(2)}$ 
into a direct sum $\oplus_{\nu = 0}^{\gaussk}  \text{ Sym}^2 M_{k + 2 \nu}$ is given. 
This leads to a basic characterization of the Spezialschar property. 
The results of this paper are directly related to the non-vanishing of certain special values of L-functions
related to the Gross-Prasad conjecture. This is illustrated by a significant example in the paper.
\end{abstract}

\ \\ \\
{\bf{Introduction}}\\
Hans Maass introduced and applied in a series of papers \cite{Ma79I},\cite{Ma79II} and \cite{Ma79III} 
the concept of a Spezialschar to prove the Saito-Kurokawa conjecture \cite{Za80}. 
Let $\mkzwei$ be the space of Siegel modular forms of degree $2$ and weight $k$. 
Let $\A$ be the set of positive semidefinite half-integral matrices of degree $2$. 
Hence $T \in \A$ can be identified with the quadratic form $T=[n,r,m]$. 
A modular form $F \in \mkzwei$ is in the Spezialschar if the Fourier coefficients 
$A(T)$ of $F$ satisfy the relation
\begin{equation}
A([n,r,m]) = \sum_{d \vert (n,r,m)} d^{k-1} A([\frac{nm}{d^2},\frac{r}{d},1])
\end{equation}
for all $\in \A$.
The space of such special forms is nowadays called the Maass Spezialschar $\maass$.
\\
\\
The purpose of this paper is twofold. First we introduce the concept of the Spezialschar 
$\mathsf{Spez}(M_k^{(2n)})$ for Siegel modular forms of even degree $2n$.
This is done in terms of the Hecke algebra ${\cal{H}}^{n}$  attached to Siegel modular forms of degree $n$.
Let us fix the embedding
\begin{eqnarray}
Sp_n \times Sp_m & \longrightarrow & Sp_{n+m} \nonumber\\
{ \tiny\vier{a}{b}{c}{d} \times \vier{\tilde{a}}{\tilde{b}}{\tilde{c}}{\tilde{d}} }& \mapsto & 
{\tiny{\vier{\tovier{a}{0}{0}{\tilde{a}}}
{\tovier{b}{0}{0}{\tilde{b}}}
{\tovier{c}{0}{0}{\tilde{c}}}
{\tovier{d}{0}{0}{\tilde{d}}}}}.
\end{eqnarray}
Let $\vert_k$ be the Petersson slash operator and let $\widetilde{T}$ be the normalized 
Heckeoperator $ T \in {\cal{H}}^n$ (see (\ref{normalize})). 
Let $\Join_T = (\widetilde{T} \times 1_{2n})- (1_{2n} \times \widetilde{T})$ and
\begin{equation}
\mathsf{Spez}\left(M_k^{(2n)}\right) := \left\{ F \in M_k^{(2n)} {\big{\vert}} \, 
F|_k \Join_T = 0 \text{ for all } T \in {\cal{H}}^n \right\}.
\end{equation}

Then we have 
\begin{Satz}
The Spezialschar introduced in this paper is the Maass Spezialschar in the case of the Siegel three-fold.
\begin{equation}
\fbox{ \text{Spez}$( M_k^{(2)}) = \maass $}.
\end{equation}
\end{Satz}
The second topic of this paper is the characterization of the space of Siegel modular forms of 
degree two and the corresponding Spezialschar in terms of Taylor coefficients and 
certain differential operators:
\begin{equation}
	\cd_{k,\underline{2\nu}}: \, \mkzwei \longrightarrow 
	M_{k+ 2 \nu}^{\mathsf{Sym}},
	\end{equation}
here $\nu \in \N_0$ and $M_{k+ 2 \nu}^{\mathsf{Sym}} = \mathsf{Sym}^2 (M_{k + 2 \nu})$.
Before we summarize the main results we give an example which also serves as an application.
Let $F_1, F_2,F_3$ be a Hecke eigenbasis of the space of Siegel cusp forms 
$S_{20}^{(2)}$ of weight $20$ and degree $2$. Let $F_1$ and $F_2$ generate the Maass Spezialschar. 
Let $f_1$ and $f_2$ be the normalized Hecke eigenbasis of $S_{24}^{(1)}$. Then we have:
\begin{equation}
\diff{20}{4}F_j = \alpha_j \, f_1 \otimes f_1 + \beta_j \, (f_1 \otimes f_2 + f_2 \otimes f_1) + \gamma_j \, f_2 \otimes f_2.
\end{equation}
It it conjectured by Gross and Prasad \cite{G-P92} that the coefficients $\alpha_j,\beta_j,\gamma_j$ 
are related to special values of certain automorphic L-functions. Recently the Gross-Prasad 
conjecture has been proven by Ikeda \cite{Ike05} for the Maass Spezialschar and $\nu=0$. 
Moreover we show in this paper that the vanishing at such special values has interesting consequences.
We have $F_j \in S_{20}^{\mathsf{Maass}}$ if and only if the special value $\beta_j$ is zero.
More generally:
		\begin{Satz}
		Let $k \in \N_0$ be even. Then we have the embedding
		\begin{equation}
		\mathbb{D}_k = \oplus_{\nu = 0}^{\gaussk}
		\cd_{k, \underline{ k + 2 \nu}} : \mkzwei \longrightarrow
		M_k^{\mathsf{Sym}} \oplus S_{k+2}^{\mathsf{Sym}} 
		\oplus \ldots \oplus
		S_{k+2 \gaussk}^{\mathsf{Sym}}.
		\end{equation}
		%
For $F\in S_k^{(2)}$ we have $\diff{k}{0}F \in \sksym$. 
\end{Satz}


Surprisingly the Maass Spezialschar property can be recovered in $M_k^{\mathsf{Sym}} \oplus S_{k+2}^{\mathsf{Sym}} 
		\oplus \ldots \oplus
		S_{k+2 \gaussk}^{\mathsf{Sym}}$ in the following transparent way. 
Let $\left( f_j \right)$ be the 
normalized Hecke eigenbasis of $M_k$. Let us define 
the diagonal subspaces
$M_k^D = \{ \sum_j \alpha_j f_j \otimes f_j \in M_k^{\mathsf{Sym}}\}$ and
$S_k^D = S_k^{\mathsf{Sym}} \cap M_k^D$. Then we can state
			\begin{Satz}
			Let $F \in \mkzwei$. Then we have
			\begin{equation}
			F \in \maass     \Longleftrightarrow   \mathbb{D}_k F \in 
			M_k^D \oplus S_{k+2}^D 
			\oplus \ldots \oplus S_{k + 2 \gaussk}^D
			\end{equation}
			and similarly
			\begin{equation}
			F \in  \smaass   \Longleftrightarrow   
			\mathbb{D}_k F \in 
			\oplus_{\nu = 0}^{\gaussk} S_{k + 2 \nu}^D.
			\end{equation}
			\end{Satz}
These two theorems give a transparent explanation of our example from a general point of view.
\\
\\
{\bf Acknowledgements:}\\
To be entered later.
\ \\ \\
{\bf{Notation}}\\
Let $Z \in \C^{n,n}$ and $\mathsf{tr}$ the trace of a matrix then we put 
$e\{Z\}=e^{2\pi i\, (\mbox{{\footnotesize tr}}\, Z)}$. 
For $l \in \Z$ we 
define $\pi_{l} = (2 \pi i)^l$. Let $x\in \R$ then we use
Knuth's notation $\knuth {x}$ to denote the 
greatest integer smaller or equal to $x$.
Let $\A_2$ denote the set of half-integral 
positive-semidefinite matrices. 
We parametrize the elements
$T = \left( \begin{smallmatrix}
			n & \frac{r}{2} \\ \frac{r}{2}  & m
			\end{smallmatrix}
			\right)$ 
			with $T=[n,r,m]$. The subset of positive-definit 
matrices we denote with $\A_2^{+}$.
\section{Ultraspherical Differential Operators}
Let us start with the notation of the ultraspherical polynomial $p_{k,\underline{2 \nu}}$. 
Let $k$ and $\nu$ be elements of $\N_0$. Let $a$ and $b$ be elements of a commutative ring.
Then we put
\begin{equation}\label{Def_ultraspherical}
\ultra{k}{2\nu}{a}{b} =
\sum_{\mu = 0}^{\nu} (-1)^{\mu} \frac{(2 \nu)!}{\mu! (2 \nu - 2 \mu)!} 
		\frac{(k+2 \nu -\mu -2)!}{(k + \nu -2)!} \, a^{2\nu - 2 \mu}\, b^{\mu}.
\end{equation}
If we specialize the parameters we have $\ultra {k}{0}{a}{b} = 1$ and 
$\ultra{k}{2\nu}{0}{0}= 0$ for $\nu \in \N $.
\\
\\
Let $\HH_n$ be the Siegel upper half-space of degree $n$. Let $\mathsf{M}_k^{(n)}$ 
the vector space of Siegel modular forms on $\HH_n$ with respect to the full modular group
$\Gamma_n = \mathsf{Sp}_n(\Z)$. Moreover let $\mathsf{S}_k^{(n)}$ denote the 
subspace of cusp forms. If $n=1$ we drop the index to simplify notation.
We denote the coordinates of the three-fold $\HH_2$ by
$(\tau, z, \tildetau )$ for
$\left( \begin{smallmatrix}
			\tau & z \\ z & \tildetau
			\end{smallmatrix}
			\right) \in \HH_2$ and 
			put $q= e\{\tau\}, \xi = e\{z\}$ and $\tilde{q}=e\{\tildetau\}$.
			Let $\widehat{d_k}$ be the dimension of $S_k$.
\begin{Def}
Let $k,\nu \in \N_0$ and let $k$ be even. Then we define the ultraspherical 
differential operator $\cd$ on the space of holomorphic functions $F$ on $\HH_2$ 
in the following way:
\begin{equation}\label{Def_differentialoperator}
\diff {k}{2 \nu} F \,
(\tau,\tildetau) \, = \, 
\left. \ultra{k}{2\nu}
{\frac{1}{2 \pi i} \frac{\partial}{\partial z}}
{\left(\frac{1}{2 \pi i}\right)^2 
\frac{\partial}{\partial \tau}
\frac{\partial}{\partial \tilde{\tau}}} F \right|_{z=0}(\tau, \tildetau).
\end{equation}
In the case $\nu =0$ we get the pullback $F(\tau,0,\tildetau)$ of $F$ on $\HH \times \HH$.
\end{Def}
Let $F \in \mathsf{M}_k^{(2)}$ with $T-\mathsf{th}$ Fourier coefficient 
$A^{F}(n,r,m)$ for $T =[n,r,m] \in \A_2$. Then we have
\begin{eqnarray} \label{diff_Fourier}
\diff {k}{2 \nu} F \, (\tau, \tildetau) & = & \sum_{n,m =0}^{\infty}
	A_{\underline{2\nu}}^F (n,m) q^n \tilde{q}^m \quad \mathsf{with} \\
	A_{\underline{2\nu}}^F (n,m) & = &
\sum_{ r \in \Z, r^2 \kleinergleich 4 nm}
\ultra {k}{2 \nu}{r}{nm} \, A^F(n,r,m). \nonumber
\end{eqnarray}
Let $\mksym = \mathsf{Sym}^2 M_k$ and $\sksym = \mathsf{Sym}^2 S_k$
$\sksym = \left( \sk \otimes \sk\right)^{\mathsf{Sym}}$. Let us 
further introduce a related Jacobi differential operator 
\textdiffjm{k}{2 \nu}. This is given by exchanging 
$\pi_{-1} \frac{\partial}{\partial \tildetau}$ with $m$ in the definition of 
the ultraspherical differential operator given 
in (\ref{Def_differentialoperator}). Applying the operator 
$\diffjm{k}{2 \nu}$ on Jacobiforms $\Phi \in \jkm$ of weight 
$k$ and index $m$ on $\HH \times \C$ matches with the effect 
of the operator $\cd_{2 \nu}$ introduced in \cite{E-Z85} ($ \S 3$, 
formula (2)) on $\Phi$. 
\\
\\
Since $F \in \mkzwei$ has a Fourier-Jacobi expansion of the form
\begin{equation}\label{Fourier_Jacobi}
F(\tau, z, \tildetau) =
\sum_{m=0}^{\infty}
\Phi_m^F(\tau,z) \, \tilde{q}^m, \,\,
\mathsf{with} \,\, \Phi_m^F \in \jkm
\end{equation}
it makes sense to consider $\diff {k}{2\nu}$ with 
respect to this decomposition in a Fourier-Jacobi expansion
\begin{equation}\label{Fourier_Jacobi Diffexpansion}
\diff{k}{2 \nu} = \bigoplus_{m=0}^{\infty} \diffjm {k}{2 \nu}.
\end{equation}
\begin{Lem}\label{welldefined}
Let $k,\nu \in \N_0$ and let $k$ be even. Then $\diff {k}{2\nu}$ maps
$\mkzwei$ to $\mksymnu$ and to $\sksymnu$ if $\nu \neq 0$. 
Moreover the subspace $\skzwei$ of cusp forms is always mapped to $\sksymnu$.
\end{Lem}
\begin{beweis}
We know from the work of Eichler and Zagier \cite{E-Z85} since 
$\diffjm {k}{2 \nu}
= \cd_{k,2\nu}$ that $\diffjm{k}{2 \nu} \Phi_m^{F} \in \mknu$.
Let $\nu >0$ then $\diffjm {k}{2 \nu} \Phi_m^F \in \sknu$ and for $F\in \skzwei$
we have $\diffjm {k}{2 \nu} \Phi_m^F \in \sknu$ for all $\nu \in \N_0$. 
We are now ready to act with
the ultraspherical differential operator 
with respect to its
Fourier-Jacobi expansion directly on the Fourier-Jacobi expansion of 
$F$ in a canonical way
\begin{equation}
\diff {k}{2 \nu} F \, (\tau, \tildetau) = \sum_{m=0}^{\infty}
\left(
\diffjm {k}{2 \nu} \Phi_m^F \right) (\tau) \, \, \tilde{q}^m,
\end{equation}
where all "coefficients"
$a_m^F (\tau) = \diffjm {k}{2 \nu} \Phi_m^F  (\tau) $ are modular forms. 
This shows us, that if we apply the Peterson slash operator 
$|_{k + 2 \nu} \gamma$ here $\gamma \in \Gamma$ to this 
function with respect to the variable $\tau$, the function is invariant. 
The same argument also works for the Fourier-Jacobi expansion 
with respect to $\tau$. From this we deduce that 
$\diff{k}{2 \nu} F \, (\tau, \tildetau) = \sum_{i,j} \alpha_{i,j} \, f_i(\tau)
f_j(\tildetau)$. Here $(f_i)_i$ is a basis of $\mknu$. Finally 
the cuspidal conditions in the lemma also follow from symmetry arguments.
\end{beweis}
\begin{Bem}
Let $F: \HH_2 \longrightarrow \C$ be holomorphic. 
Let $g \in Sl_2(\R)$ and let 
$J = \left( \begin{smallmatrix}
			0 & 1_2\\ -1_2  & 0
			\end{smallmatrix}
			\right)$.
Then we have:
\begin{eqnarray}
\diff{k}{2 \nu}(F\vert_k (g\times 1_2)) & = & (\diff{k}{2 \nu}F)\vert_{k+ 2 \nu} (g\times 1_2)\\
\diff{k}{2 \nu}(F\vert_k (1_2 \times g)) & = & (\diff{k}{2 \nu}F)\vert_{k+ 2 \nu} (1_2 \times g)\\
\diff{k}{2 \nu}(F\vert_k J) & = & (\diff{k}{2 \nu}F)\vert_{k+ 2 \nu} J.
\end{eqnarray}
\end{Bem}

\begin{Bem}
There are other possibilities for construction of
differential operators as used in this section (see Ibukiyama for a overview \cite{Ibu99}). 
But since the connection between our approach and the theory 
developped of Eichler and Zagier \cite{E-Z85} is so useful we 
decided to do it this way. We also wanted to introduce the 
concept of Fourier-Jacobi expansion of differential 
operators, which is interesting in its own right.
\end{Bem}

\section{Taylor Expansion Of Siegel Modular Forms}

The operators $\diff {k}{2\nu}$ can be seen at this point
as somewhat artificial. If we apply $\diffknu$ to 
Siegel modular forms $F$ we lose information. For example we know that 
$\mathsf{dim}S_{20}^{(2)} = 3$ and contains a two dimensional 
subspace of Saito-Kurokawa lifts. 
Since $\mathsf{dim} S_{20}^{\mathsf{Sym}} =1$ we obviously lose 
informations if we apply $\diff {20}{0}$. But even worse let 
$F_1$ and $F_2$ be a Hecke eigenbasis of the space of 
Saito-Kurokwa lifts and $F_3$ a 
Hecke eigenform of the orthogonal complement then we 
have $\diff {20}{0}F_j \neq 0$ for $j= 1,2,3$. 
The general case seems to be even worse, since for example
$\dimens \mkzwei \sim k^3$ and $\dimens \mksym \sim k^2$. On the other 
hand from an optimistic viewpoint we may find about 
$k$ pieces $\diffknu F$ which code all the relevant information needed to 
characterize the Siegel modular forms $F$.
\\
\\
Paul Garrett in his fundamental papers \cite{Ga84} and \cite{Ga87} 
introduced the method of calculating pullbacks of 
modular forms to study 
automorphic L-functions. We also would like to mention the work of Piatetski-Shapiro, 
Rallis and Gelbart at this point (see also \cite{GPR87}).
And recently Ichino in his 
paper: Pullbacks of Saito-Kurokawa lifts \cite{Ich05} 
extended Garrett's ideas in a brilliant way to prove the Gross-Prasad 
conjecture \cite{G-P92} for Saito-Kurokawa lifts. 
In the new language we have introduced, it is obvious to 
consider Garretts pullbacks as the $0-$th 
Taylor coefficients of $F$ around $z=0$. Hence it seems to 
be very lucrative to study also the higher Taylor 
coefficients and hopefully get some transparent link.
\\
\\
Let $k \in \N_0$ be even. Let $F \in \mkzwei$ and $\Phi \in \jkm$. 
Then we denote by
\begin{equation}\label{Taylor}
F(\tau, z, \tildetau) = \sum_{\nu = 0}^{\infty} 
\chi_{2 \nu}^F (\tau, \tildetau) \, \, z^{2\nu}
\quad \mathsf{and}\quad
\Phi (\tau,z) = \sum_{\nu = 0}^{\infty} \chi_{2 \nu}^{\Phi} (\tau) \, \, z^{2\nu}
\end{equation}
the correponding Taylor expansions with respect to $z$ around $z=0$. 
Here we already used the invariance of $F$ and $\Phi$ with 
respect to the transformation $z \mapsto (-z)$ since $k$ is even. 
Suppose $\chi_{2 \nu_0}$ is the first non-vanishing Taylor coefficient, 
then we denote $2\nu_0$ the vanishing order of the underlying form. 
If the form is identically zero we define the vanishing order to be $\infty$.
To simplify our notation we introduce normalizing factor
\begin{equation}\label{normalizing}
\gamma_{k, \underline{2\nu}} = \left( \frac{1}{2 \pi i}\right)^{2\nu} 
\frac{
(k + 2 \nu -2)! \,\, (2\nu)!
}
{(k+\frac{2\nu}{2}-2)!
}.
\end{equation}
Further we put
\begin{equation}\label{muliti_differential}
\chi_{2\nu}^{\mu,\mu} =
\frac{\partial^{2\mu}}
				{\partial \tau^{\mu}\partial \tilde{\tau}^{\mu}} \chi_{2\nu }^F
\,\, \,\,\mathsf{and}\,\,\,\,
\xi_{2\nu}^{\mu,\mu} =
\left(\gamma_{k, \underline{2\nu}}\right)^{-1}
\frac{\partial^{2\mu}}
				{\partial \tau^{\mu}\partial \tilde{\tau}^{\mu}} 
\cd_{k, \underline{2\nu}}F.				
\end{equation}
Then a straightforward calculation leads to the following useful formula.
\begin{Lem}
Let $k, \nu \in \N_0$ and let $k$ be even. Let $F \in \mkzwei$. Then we have
\begin{equation}\label{link_Diff_Taylor}
\left( \cd_{k, \underline{2 \nu}}F \right) (\tau, \tilde{\tau})
=
\gamma_{k,\underline{2 \nu}} \sum_{\mu =0}^{\nu} (-1)^{\mu}
\frac{(k+2 \nu - \mu -2)!}{(k + 2 \nu -2)!\mu!}
\left(\frac{
\partial^{2\mu} \chi_{ 2\nu - 2\mu}^F}{\partial \tau^{\mu}\,\, 
\partial \tilde{\tau}^{\mu}}	\right).
\end{equation}
A similiar formula is valid for Jacobiforms with 
normalizing factor $\gamma_{k, \underline{2\nu}}^{\mathsf{J},m} = 
\gamma_{k, \underline{2\nu}}$.
\end{Lem}
\begin{Cor}
Let $2 \nu_0$ be the vanishing order of $F \in \mkzwei$. 
Then we have $\diffknu F =0$ for $\nu < \nu_0$ and
\begin{equation}\label{identiy_vanishingorder}
\diff{k}{2\nu_0} F \, (\tau, \tildetau) = 
\gamma_{k, \underline{2\nu}}\,\,
\chi_{2 \nu_0}^F (\tau, \tildetau) \in \mksymnunull
\backslash \{0\} .
\end{equation}
\end{Cor}
Similiary we have for $\Phi \in \jkm$ with vanishing order 
$2 \nu_0$ the properties $\diffjm{k}{2 \nu}\Phi = 0$ for 
$\nu < \nu_0$ and $\diffjm{k}{2 \nu_0}\Phi  = 
\gamma_{k, \underline{2\nu}}\,\,
\chi_{2 \nu_0}^{\Phi} \in M_{k + 2 \nu_0}$. 
\\
\\
\underline{{\bf EXAMPLE:}}
It is well known that $\dimens\,  S_{10}^{(2)} =1$. 
Let $\Phi = \Phi_{10} \in S_{10}^{(2)} $
be normalized in such a way that $A^{\Phi}(1,1,1) =1$. 
Then it follows from $\diff {10}{0} \Phi = 0$ that $A^{\Phi}(1,0,1) = -2$ 
since $\dimens \,S_{10}^{\mathsf{Sym}} = 0$. 
Then $\Phi$ has the Taylor expansion
\begin{equation}\label{example_ten}
\Phi_{10}(\tau, z, \tildetau) =
\frac{3}{5} \pi_2 \Delta(\tau) \Delta(\tildetau) \,\, z^2 +
\Delta{\mathsf{'}}(\tau) \Delta{'}(\tildetau) \,\, z^4 + O(z^6).
\end{equation}
\\
We can also express the Taylor coefficients $\chi_{2 \nu}^F$ 
in terms of the modular forms $\diff{k}{2\nu}F$. This can be done 
by inverting the formula (\ref{link_Diff_Taylor}). Finally we get 
				\begin{equation}\label{Essenz}
				\fbox{$\chi_{2\nu} = 
				\sum_{\mu = 0}^{\nu}
				\frac{(k+2 \nu - 2\mu -1)!}
				{(k+2 \nu -\mu-1)! \mu!}
				\,\,\xi_{2\nu - 2 \mu}^{\mu,\mu} .
				$}
				\end{equation}
			
Before we state our first main result about the entropy of the 
family $\diff {k}{0}F$, $\diff {k}{2}F$, $\diff {k}{4}F \ldots $ 
we introduce some further notation. 
					\begin{eqnarray}\label{directsum}
					\mathbb{W}_k & = & M_k^{\mathsf{Sym}} \oplus 
					\bigoplus_{j=1}^{\gauss{\frac{k}{10}}} S_{k +2j}^{\mathsf{Sym}}
					\qquad \mathsf{and} \\
					\mathbb{W}_k^{\mathsf{cusp}} & = & S_k^{\mathsf{Sym}} \oplus 
					\bigoplus_{j=1}^{\gauss{\frac{k}{10}}} S_{k +2j}^{\mathsf{Sym}}.
					\end{eqnarray}
These spaces will be the target of our next consideration. 
More precisely we define a linear map from the space of Siegel modular forms of degree $2$ into these spaces with 
remarkable properties.

\begin{Satz}\label{embedding}
Let $k \in \N_0$ be even. Then we have the linear embedding
\begin{equation}
\mathbb{D}_k : 
\begin{cases}
M_k^{(2)} \hookrightarrow \mathbb{W}_k \\
F \mapsto
\bigoplus_{\nu =0}^{\gauss{\frac{k}{10}}} 
\cd_{k, \underline{2 \nu}}F .
\end{cases}
\end{equation}
Since $\diff {k}{0} \skzwei$ is cuspidal we have 
the embedding of $\skzwei$ into $\mathbb{W}_k^{\mathsf{cusp}} $.
\end{Satz}
\begin{Bem}
It can be deduced from \cite{Hei06} that 
$\diff{k}{0} \oplus \diff {k}{2}$ is surjective.
Hence for $k < 20$ we have:
\begin{itemize}
\item
$\mkzwei$ is isomorphic to $\mk$ for $k< 10$ and
\item
$\mkzwei$ is isomorphic to $\mksym \oplus \sksym$ for 
$10 \kleinergleich  k < 20$ and $\skzwei \simeq \sk \oplus S_{k+2}$.
\end{itemize}
\end{Bem}

\begin{beweis}
First of all we recall that we have already shown that 
$\diff {k}{0} \mkzwei \subseteq \mksym$ and 
$\diff {k}{2\nu} \mkzwei \subseteq \sksymnu$ for $\nu > 0$.
Let $F\in \mkzwei$ and suppose that $\mathbb{D}_k F$ is identically zero. 
Then it would follow from our inversion formula (\ref{Essenz}) that
\begin{equation}\label{higher}
F(\tau, z, \tildetau) = \sum_{\nu =\gauss{\frac{k}{10}}+1}^{\infty}
\chi_{2 \nu}^F(\tau, \tildetau) \,\, z^{2\nu}.
\end{equation}
For such $F$ the general theory of Siegel modular forms of degree $2$ 
says that the special function $\Phi_{10} \in \skzwei$, 
which we already studied in one of our examples, divides $F$ 
in the $\C-$algebra of modular forms. And this is 
fullfilled at least with a power of $\gauss{\frac{k}{10}}+1= t_k >0$. 
Hence there exists a Siegel modular form $G$ of weight $k - 10\, t_k$. 
But since this weight is negative and non-trival Siegel modular 
forms of negative weight do not exist the form $G$ has to 
be identically zero. Hence we have shown that if 
$\mathbb{D}_k F \equiv 0$ then $F \equiv 0$. And this proves the 
statement of the theorem.
\end{beweis}
\begin{Bem}
The number $\gauss{\frac{k}{10}}$ in the Theorem is optimal. This follows directly from properties of $\Phi_{10}$.
\end{Bem}
\begin{Bem}
Let $E_{k}^{2,1}(f)$ be a Klingen Eisenstein series attached to $f \in S_k$. 
Let $E_k$ denote an elliptic Eisenstein series of weight $k$.
Then it can be deduced from \cite{Ga87} that $\diff{k}{0}E_k^{2,1}(f) = 
f \otimes E_k + E_k \otimes f $ mod $\sksym$. 
\end{Bem}
\begin{Bem}
It would be interesting to have a different proof of the Theorem independent of the special properties of $\Phi_{10}$. 
\end{Bem}
\begin{Bem}
The asymptotic limit of the dimension of the quotient of 
$\mathbb{W}_k  \, / \, \mkzwei$ is equal to $\frac{91}{25}$. Let us put $d_k = \text{ dim }M_k$.
\begin{itemize}
\item
The dimension of the target space $\mathbb{W}_k$:
\begin{eqnarray*}
\vspace*{-3cm}
\mathsf{dim} \, \mathbb{W}_k  & \sim & 
\frac{1}{288} \,\, \int_{0}^{\frac{k}{10} } (k + 2x)^2 \, dx\\
& \sim & \frac{1}{288} \,\,\frac{1}{2 \cdot 3} \,\,\frac{91}{5^3} \,\, k^3
\end{eqnarray*}
\item
The asymptotic dimension formula of $\mkzwei$ is given by
\begin{eqnarray*}
\ \vspace*{-3cm}
\phantom{xxxx}\mathsf{dim}\, \mkzwei & \sim &  \frac{1}{288} \,\,\frac{1}{2 \cdot 3 \cdot 5}  \,\, k^3 \text{ (see \cite{Ma79I}, Introduction)}.
\end{eqnarray*}
\end{itemize}
\end{Bem}

\section{The Spezialschar}
In this section we first recall some basic facts on the Maass Spezialschar \cite{Za80}. 
Then we determine the image of the 
Spezialschar in the space $\mathbb{W}_k$ for all even weights $k$. 
Then finally we introduce a 
Spezialschar as a certain subspace of the space of Siegel 
modular forms of degree $2n$ and weight $k$. Then we show that in the case $n=1$ 
this Spezialschar coincides with the Maass Spezialschar .

\subsection{Basics of the Maass Spezialschar}
Let $\jkm$ be the space of Jacobi forms of weight $k$ and 
index $m$. We denote the subspace of cusp forms with $\jkmcusp$. Let $\mid_{k,m}$ 
the slash operator for Jacobi forms and $V_l$ ($l \in \N_0$) be the operator, which maps
$\jkm$ to ${{J}}_{k,ml}$. More precisely, 
let $\Phi (\tau,z) = \sum c(n,r)\, q^n \xi^r \in \jkm$.
Then $(\Phi \mid_{k,m} V_l)(\tau,z) = 
\sum c^{*}(n,r)\, q^n \xi^r$ with
\begin{equation}
c^{*}
(n,r) = \sum_{a \mid (n,r,l)} a^{k-1} \, c(\frac{nl}{a^2},\frac{r}{a})\quad 
\mathsf{ for} \,\, l\in \N
\end{equation}
and for $l=0$, we have 
$c^{*}(0,0) = c(0,0) \,\, \left( \frac{-2k}{B_{2k}}\right)$
and for $l=0$ and $n > 0$ we have $c^{*}(n,r)= c(0,0) \,\, \sigma_{k-1}(n)$. 
This includes the theory of Eisenstein series in a nice way \cite{E-Z85}.


\begin{Def}
The 
lifting
$\cal{V}$ is given by the linear map
		\begin{equation}
		{\cal{V}} : 
					\begin{cases}
					{{J}}_{k,1} \longrightarrow M_k^{(2)} \\
					\Phi \mapsto \sum_{l=0}^{\infty} \left(\Phi \mid_{k,1} 
					V_l \right) \widetilde{q}^l.
					\end{cases} \label{SaitoKurokawa}
		\end{equation}
The image of this lifting is the Maass 
Spezialschar $M_k^{\mathsf{Maass}}$
of weight $k$. 
The subspace of cusp forms we denote with
$S_k^{\mathsf{Maass}}$.
\end{Def}

\begin{Bem}\
\begin{itemize}
\item
The lifting is invariant by 
the Klingen parabolic of $ Sp_2(\Z)$. Since the Fourier coefficients satisfy
$A(n,r,m) = A(m,r,n)$ the map ${\cal{V}}$ is well-defined.
\item
If we restrict the Saito-Kurokawa lifting to Jacobi cusp forms we get Siegel cusp forms.
\item
Let $\Phi \in \jkm$ and $l, \mu \in \N_0$. Then we have
			\begin{equation}
			\cd_{k, \underline{ {2 \mu}}}^{J,ml} \left( \Phi\mid_{k,m}V_l\right)
			=
			\left(
			\cd_{k, \underline{ {2 \mu}}}^{J,m} \Phi
			\right)
			\mid_{k}T_l.
			\label{interchange}
			\end{equation}
			Here $T_l$ is the Hecke operator on the space of elliptic modular forms.
\item
Let $F \in M_k^{\mathsf{Maass}}$ be the lift of $\Phi \in J_{k,1}$. 
Then $F$ is a Hecke eigenform if and only if $\Phi$ is a Hecke-Jacobi eigenform.
\end{itemize}
\end{Bem}
From this consideration we conclude \cite{E-Z85}:
\begin{Prop}
Let $F \in M_k^{(2)}$ be a Siegel modular form. Then the following properties are
aquivalent
\begin{itemize}
\item
\underline{ARITHMETIC} Let $A(n,r,m)$ denote the Fourier coefficients of $F$ then
\begin{equation}
A(n,r,m) = \sum_{d\mid (n,r,m)} d^{k-1} A(\frac{n\, m}{d^2}, \frac{r}{d},1)
\end{equation}
\item
\underline{LIFTING}
Let $\Phi_1^F$ be the first Fourier-Jacobi coefficient of $F$. 
Then all other Fourier-Jacobi coefficients satisfy the identity
\begin{equation}
\Phi_m^F = \Phi_1^F \mid_{k,1} V_m \,\,.
\end{equation}
\end{itemize}
Let $F\in \skzwei$ be a Hecke eigenform. 
Then $F$ is a Saito-Kurokawa lift if and only if the spinor 
L-function $Z(F,s)$ of degree $4$ has a pole (\cite{Ev80}).
\end{Prop}
\subsection{The Diagonal of $\mathbb{W}_k$}
Let $(f_j)$ be the normalized Hecke eigenbasis of $M_k$. With this 
notation we introduce the diagonal space
\begin{equation}
\mkd = \{ \sum_j \alpha_j \, f_j \otimes f_j \in \mksym \}
\end{equation}
and the corresponding cuspidal subspace $\skd$. Now we are ready to 
distinguish the Maass Spezialschar in the vector 
spaces $\mathbb{W}_k$ and $\mathbb{W}_k^{\mathsf{cusp}}$.

\begin{Satz}\label{diagonal}
Let $k$ be a natural even number. Let $F$ be a Siegel modular form 
of degree two and weight $k$. Then we have
		\begin{equation}
		F \in \maass \Longleftrightarrow \D_k F \in \mkd \oplus 
		{S_{k+ 2 }^{D}}  \oplus
		\ldots
		\oplus
		S_{k+ 2 \gauss{\frac{k}{10}}}^{D}.
		\end{equation}
Let $F$ be a cuspform. Then we have
		\begin{equation}
		F \in \smaass \Longleftrightarrow \D_k F \in 
		\oplus_{\nu =0}^{ \gauss{\frac{k}{10}}} S_{k + 2 \nu}^{D}.
		\end{equation}
\end{Satz}
\ \\ 
\begin{Bem}\
\begin{itemize}
\item
The Theorem \ref{diagonal} describes a link between Siegel modular forms and elliptic Hecke eigenforms.
\item
Let $F\in M_{20}^{(2)}$ and let $(f_j)$ be a Hecke eigenbasis of $S_{24}$. Then $ F\in \maass$ if and only if
\begin{equation}
\diff{20}{4} F = \alpha_0 E_{24} \otimes E_{24} + \alpha f_1 \otimes f_1 + \gamma f_2 \otimes f_2
\end{equation}
here $\alpha_0, \alpha, \gamma \in \C$.
\end{itemize}

\end{Bem}
\begin{beweis}
We first show that if $F$ is in the Maass Spezialschar then 
$\diff{k}{2 \nu} F$ is an element of the diagonal space.
Let $\nu \in \N_0$ and $\Phi_1^F$ be the first Fourier-Jacobi 
coefficient of $F$. Then we have
\begin{equation}
				\left( 
					\cd_{k, \underline{2 \nu}} \left( {\cal{V}} \Phi\right) \right) 
				(\tau, \tilde{\tau}) =
				\sum_{l=0}^{\infty} 
				\left( \cd_{k, \underline{2 \nu}}^{J,l}
				\left( \Phi \vert_{k,1} V_l \right)\right) (\tau) \,\, \widetilde{q}^l.
				\end{equation}
Here we applied the Fourier-Jacobi expansion of the differential operator 
$\cd_{k, \underline{2 \nu}} $ acting on Siegel modular forms. Then we used the formula
(\ref{interchange})
to interchange the operators
$\cd_{k, \underline{2 \nu}}^{J,l} $ and $V_l$ to get
\begin{equation}
\left(\cd_{k, \underline{2 \nu}} F \right) (\tau, \tilde{\tau}) =
\sum_{l=0}^{\infty}
\left(
\cd_{k, \underline{2 \nu}}^{J,l} \Phi
\right)\vert_{k} T_l \,\,
\widetilde{q}^l.
\end{equation}
Now let
$\left( f_j^{\underline{k + 2 \nu }}\right)_{j=1}^{\hutvontknu}$ 
be a normalized Hecke eigenbasis of $S_{k + 2 \nu}$. Let 
$1 \le j_1,j_2 \le \hat{d}_{k + 2 \nu}$. Then we have
\begin{equation} \label{reprod}
\llk
\left(\cd_{k, \underline{2 \nu}} F\right)
(\,{*}\, , \tilde{\tau}), 
f_{j_1}^{\underline{k + 2 \nu }}
\rrk
=
\llk
\left(\cd_{k, \underline{2 \nu}}^{J,l} \Phi\right) ,
f_{j_1}^{\underline{k + 2 \nu }}
\rrk
\,\,
f_{j_1}^{\underline{k + 2 \nu }},
\end{equation}
which leads to the desired result
\begin{equation}
\llk
\left(\cd_{k, \underline{2 \nu}} F\right),
f_{j_1}^{\underline{k + 2 \nu }}
\otimes
f_{j_2}^{\underline{k + 2 \nu }}\rrk
= 
0
\,\,
\mathsf{for}\,\,
j_1 \neq j_2.
\end{equation}
It remains to look at the Eisenstein part if $\nu =0$. 
Since the space of Eisenstein series has the basis $E_k$ and
is orthogonal to the functions given in (\ref{reprod}) we have proven that the 
Spezialschar property of $F$ implies that $\mathbb{D}_k F \in
\mathbb{W}_k ^{\mathsf{D}}$.\\
\\
Now let us assume that $F \notin M_k^{\mathsf{Maass}}$. 
Then we show that $\mathbb{D}_k F \notin
\mathbb{W}_k ^{\mathsf{D}}$. Since the map 
\begin{equation}
(\cd_{k, \underline{0}} \oplus
\cd_{k, \underline{2 }} ): \maass \longrightarrow M_k^D \oplus S_{k+2}^D
\end{equation}
is an isomorphism, we can assume that
$\left( \cd_{k, \underline{0}} \oplus \cd_{k, \underline{2 }} \right)
(F)$ projected on 
$M_k^D \oplus S_{k+2}^D$ is identically zero. Altering $F$ by an element 
of the Maass Spezialschar does not change the property we have to prove. 
If $\cd_{k, \underline{0}} F \notin M_k^D$ or $\cd_{k, \underline{2}} F
\notin S_{k +2}^D$ we are done otherwise we can assume that
\begin{equation}
\left(\cd_{k, \underline{0}} \oplus
\cd_{k, \underline{2 }}\right)(F) \equiv 0.
\end{equation}
Then we have the $\mathsf{order}F = 2 \nu_0 \grossergleich 4$ and $k \grossergleich 20$, since
$F \notin M_k^{\mathsf{Maass}}$. Let
\begin{equation}
F \vier{\tau}{z}{z}{\tilde{\tau}} =
\sum_{\nu = \nu_0}^{\infty}
\chi_{2 \nu_0}^F (\tau, \tilde{\tau})\,\,  z^{2 \nu}
\end{equation}
be the Taylor expansion of $F$ with 
$\chi_{2 \nu_0}^F(\tau, \tilde{\tau}) \in S_{k + 2 \nu_0}$ not identically zero.
Let $\Phi_{10} \in S_{10}^{(2)}$ be the Siegel cusp form (\ref{example_ten}) 
of weight $10$ and degree $2$. It has the properties that
$\chi_0^{\Phi_{10}} \equiv 0$ and $\chi_2^{\Phi_{10}}(\tau, \tildetau) =
c \,\,
\Delta(\tau) \,\Delta(\tilde{\tau})$ with $c\neq 0$. Since 
$\mathsf{order}F = 2 \nu_0$ we also have
		\begin{equation}
		\Phi_{10}^{\nu_0} \parallel F.
		\end{equation}
This means that there exists a $G \in S_{k-10\nu_0}$ such that $\chi_0^G$ is non-trivial and
		\begin{equation}
		F = \left( 
				\Phi_{10}
								\right)^{\nu_0} \,\, G.
		\end{equation}
Hence we have for the first nontrivial Taylor coefficient of $F$ the formula
\begin{eqnarray}
\chi_{2 \nu_0}^F(\tau, \tilde{\tau})  &=& \left( \chi_2^{\Phi_{10}} 
(\tau, \tilde{\tau}) \right)^{\nu_0} \,\, \chi_0^G (\tau, \tilde{\tau}) \\
&=&
c^{\nu_0} \,\, \Delta(\tau)^{\nu_0}  \, \Delta(\tilde{\tau})^{\nu_0} 
\chi_0^G (\tau, \tilde{\tau}).
\end{eqnarray}
And the coefficient $a_1(\tilde{\tau})$ of $q$ is identically zero. 
Now let us assume for a moment that $\chi_{2 \nu_0}^F \in S_{k+ 2 \nu_0}^D$. Then we have
\begin{equation}
\chi_{2 \nu_0}^F (\tau, \tilde{\tau}) = \sum_{l=1}^{\hat{d}_{k + 2 \nu}} \alpha_l \,\, 
f_l^{\underline{k + 2 \nu_0}}(\tau) \,\,
f_l^{\underline{k + 2 \nu_0}}(\tilde{\tau}) 
\end{equation}
and the coefficient of $q$ is given by $
\sum_{l=1}^{\hat{d}_{k+2\nu_0}} \alpha_l \,\, 
f_l^{\underline{k + 2 \nu_0}}(\tilde{\tau})$. Since 
$\left( f_l^{\underline{k + 2\nu_0}}\right)_{l=1}^{\hat{d}_{k+2 \nu_0}}$ 
is a basis we have $\alpha_1 = \ldots = \alpha_{\hat{d}_{k+2 \nu_0}} = 0$.
But since we assumed that $\mathsf{order}F = 2 \nu_0$ we have 
a reductio ad absurdum. Hence we have shown that 
$\chi_{2 \nu_0}^F \notin S_{k+2\nu_0}$ which proves our theorem.
\end{beweis}
\begin{Cor}
Klingen Eisenstein series are not in the Maass Spezialschar.
\end{Cor}
\begin{Bem}\label{general}
Let $k$ be a natural even number. Let $F$ be a Siegel modular form 
of degree two and weight $k$. Then we have
		\begin{equation}
		F \in \maass \Longleftrightarrow \diff{k}{2 \nu} F \in \mkdnu \text{ for all } \nu \in \N_0.
		\end{equation}
\end{Bem}
\subsection{The Spezialschar}
Let $G^{+}Sp_n(\Q)$ be the rational symplectic group 
with positive similitude $\mu$. 
In the sense of Shimura we attach to 
Hecke pairs the corresponding Hecke algebras
\begin{eqnarray}
{\cal{H}}^n &=& \left( \Gamma_n , G^{+}Sp_n(\Q)\right)\\
{\cal{H}}_0^n &=& \left( \Gamma_n , Sp_n(\Q)\right).
\end{eqnarray}	
We also would like to mention that in the setting of 
elliptic modular forms the classical Hecke operator $T(p)$
can be normalized such that it is an element of the full Hecke algebra ${\cal{H}}^1$, 
but not of the even one ${\cal{H}}_0^1$.
Let $g \in G^{+}Sp_n(\Q)$ with similitude $\mu (g)$. Then we put
\begin{equation}\label{normalize}
\widetilde{g} = \mu(g)^{-\frac{1}{2}} \, g
\end{equation}
to obtain an element of $Sp_n(\R)$. We further extend this to ${\cal{H}}^n$.
\begin{Def}
Let $T \in {\cal{H}}^n $. Then we define
\begin{equation}
\Join_T = (\widetilde{T} \times 1_{2n})- (1_{2n} \times \widetilde{T}).
\end{equation}
Here $\times$ is the standard embedding of $(Sp_n, Sp_n)$ into $Sp_{2n}$.
\end{Def}
Now we study the action $\vert_k \Join_T $ on the space of 
modular forms of degree $2n$ for all $T \in {\cal{H}}^n $ or 
$T \in {\cal{H}}_0^n $. The first thing we would like to mention is that 
for $F \in M_k^{(2n)}$ the function
$F\vert_k \Join_T $ is in general not an element of 
$M_k^{(2n)}$ anymore. Anyway at the moment 
we are much more interested in the properties 
of the kernel of a certain map related to this action.
In particular in the case $n=1$ we get a new description of the Maass Spezialschar.

\begin{Def}
Let $n$ and $k$ be natural numbers. Let $M_k^{(2n)}$ be the space of Siegel 
modular forms of degree $2n$ and weight $k$. Then we introduce the Spezialschar corresponding
to the Hecke algebras ${\cal{H}}^n$ and  ${\cal{H}}_0^n$.

\begin{eqnarray}
\mathsf{Spez}\left(M_k^{(2n)}\right) \!\!\!& = & \!\!\!\left\{ F \in M_k^{(2n)} {\big{\vert}} \, 
F|_k \Join_T = 0 \text{ for all } T \in {\cal{H}}^n \right\}\\
\mathsf{Spez}_0\left(M_k^{(2n)}\right)\!\!\! & = & \!\!\!\left\{ F \in M_k^{(2n)} {\big{\vert}} \, 
F|_k \Join_T = 0 \text{ for all } T \in {\cal{H}}_0^n \right\}.
\end{eqnarray}
Moreover $\mathsf{Spez}\left(S_k^{(2n)}\right) $ and $\mathsf{Spez}_0\left(S_k^{(2n)}\right) $ 
are the cuspidal part of the correponding Spezialschar.
\end{Def}
It is obvious that these subspaces of $M_k^{(2n)}$ are candidates for 
finding spaces of modular forms with distinguished Fourier coefficients. 
Further it turns that these spaces are related to the Maass Spezialschar and the Ikeda lift \cite{Ike01}. 
More precisely in the first interesting case we have:

\begin{Satz}\label{equality}
The Spezialschar $\mathsf{Spez}\left(M_k^{(2)}\right)$ is equal to the Spezialschar of Maass.
\end{Satz}
\begin{beweis}
Let $F\in M_k^{(2)}$. Then we have $F \in \maass$ if and only if 
$\diff{k}{2 \nu} F \in \mkdnu$ for all $\nu \in \N_0$. This follows from Remark \ref{general}. 
On the other side the property $\diff{k}{2 \nu} F \in \mkdnu$ is equivalent to the identity
\begin{equation}
\left( \diff{k}{2 \nu} F \right)\vert_{k+ 2 \nu} \Join_T = 0 \text{ for all } T \in {\cal{H}}.
\end{equation}
This follows from the fact that the Hecke operators are self adjoint and that the 
space of elliptic modular forms has multiplicity one.
To make the operator well-defined we used the embedding 
$\HH \times \HH$ into the diagonal of $\HH_2$. 
We can now interchange the differential operators $\diff{k}{2\nu}$ and the Petersson slash operator $\vert_{*}$.
This leads to
\begin{equation}
\diff{k}{2 \nu} F \in \mkdnu  \Longleftrightarrow  
\diff{k}{2 \nu} \left( F
\vert_{k} \Join_T \right)= 0.
\end{equation}
So finally it remains to show that if $\diff{k}{2 \nu} \left( F
\vert_{k} \!\! \Join_T \right)=0$ for all $\nu \in \N_0$ then it follows 
$F\vert_k \!\! \Join_T =0$. 
By looking at the Taylor expansion of the function 
$F\vert_k \Join_T \left(\ovier{\tau}{z}{z}{\tildetau}\right)$ 
with respect to $z$ around $0$ we get with the same 
argument as given in the proof of Theorem \ref{embedding} the desired result.
\end{beweis}
\section{Maass relations revised}
We introduced two Hecke algebras ${\cal{H}}$ and ${\cal{H}}_0$ related to elliptic modular forms.
For the correponding Spezialschar $\mathsf{Spez}(M_k^{(2)})$ and $\mathsf{Spez}_0(M_k^{(2)})$ we obtain:

\begin{Satz}
Let $k$  be an even natural number. Then the even Spezialschar $\mathsf{Spez}_0(M_k^{(2)})$ related to the Hecke algebra 
${\cal{H}}_0$ which is locally generated by $T(p^2)$ is equal to the 
Spezialschar $\mathsf{Spez}(M_k^{(2)})$ related to the Hecke algebra ${\cal{H}}$ which is locally generated by $T(p)$.
\begin{equation}
\mathsf{Spez}_0(M_k^{(2)}) = \mathsf{Spez}(M_k^{(2)}).
\end{equation}
\end{Satz}
\begin{beweis}
Let $F\in M_k^{(2)}$. 
We proceed as follows. In the proof of Theorem \ref{equality} it has been shown that
\begin{equation}
F \in \mathsf{Spez}\left(M_k^{(2)}\right)
\Longleftrightarrow  
\left( \diff{k}{2 \nu} F \right)\vert_{k+ 2 \nu} \Join_T = 0 \text{ for all } T \in {\cal{H}} 
\text{ and }
\nu \in \N_0.
\end{equation}
Now we show that
\begin{equation}\label{left}
\left( \diff{k}{2 \nu} F \right)\vert_{k+ 2 \nu} \Join_{T(p)}= 0 
\Longleftrightarrow  
\left( \diff{k}{2 \nu} F \right)\vert_{k+ 2 \nu} \Join_{T(p^2)}= 0 
\end{equation}
for all $\nu \in \N_0$ and prime numbers $p$. This would finish the proof since
\begin{equation}
F \in \mathsf{Spez}_0\left(M_k^{(2)}\right)
\Longleftrightarrow  
\left( \diff{k}{2 \nu} F \right)\vert_{k+ 2 \nu} \Join_T = 0 \text{ for all } T \in {\cal{H}}_0
\text{ and }
\nu \in \N_0.
\end{equation}
(this can also be obtained by following the procedure of the proof of Theorem \ref{equality}).\\ \\
To verify the equation (\ref{left}) we show that to being an element of the 
kernel of the operator $\vert \Join_{T(p^2)}$ implies already to be an element of the kernel of $\vert \Join_{T(p)}$. \\ \\
To see this we give a more general proof. Let $\phi \in \mksym$ and let $\phi \vert_k \Join_{T(p^2)} = 0$. Let $(f_j)$ be a
normalized Hecke eigenbasis of $M_k$. Then we have
\begin{equation}
\phi = 
\sum_{i,j} \alpha_{i,j} \, f_i \otimes f_j
\end{equation}
Let us assume that there exists a $\alpha_{i_0,j_0} \neq 0$ with $i_0 \neq j_0$. 
Let us denote $\lambda_l(p^2)$ to be the eigenvalue of $f_l$ with respect to the Hecke operator $T(p^2)$. Then we have
\begin{equation}
0 = \phi \vert_k \Join_{T(p^2)} = \sum_{i,j} \alpha_{i,j} ( \lambda_i(p^2) - \lambda_j(p^2)        )\, f_i \otimes f_j.
\end{equation}
From this follows that $\lambda_{i_0}(p^2) =\lambda_{j_0}(p^2) $ for all prime numbers $p$.  
It is easy to see at this point that then $f_{i_0}$ and $f_{j_0}$ have to be cusp forms. 
In the setting of cusp forms we can apply a result on multiplicity one for 
$SL_2$ of D. Ramakrishnan \cite{Ra00}(section 4.1) and other people to obtain $f_{i_0} = f_{j_0}$. 
Since this is a contradiction we have
$\phi \in \mkd$. In other words we have $\phi \vert_k \Join_{T(p)} = 0$.
\end{beweis}

\
\\
\
\begin{center}
\mbox{----------------------}
\end{center}
\
\\
\
\begin{minipage}{6.5cm}
{\it Bernhard Heim}\\
MPI Mathematik\\
Vivatsgasse 7\\
53111 Bonn, Germany\\
heim\textcircled{a}mpim-bonn.mpg.de
\end{minipage}
\end{document}